\documentclass{ieeeoj}
\usepackage{microtype}
\usepackage{cite}
\usepackage{amsmath,amssymb,amsfonts}
\usepackage{algorithmic}
\usepackage{graphicx,color}
\usepackage{textcomp}
\usepackage{hyperref}
\def\BibTeX{{\rm B\kern-.05em{\sc i\kern-.025em b}\kern-.08em
    T\kern-.1667em\lower.7ex\hbox{E}\kern-.125emX}}
\AtBeginDocument{\definecolor{green}{cmyk}{1,0,1,0.5}}
\def\OJlogo{\vspace{-14pt}\includegraphics[height=28pt]{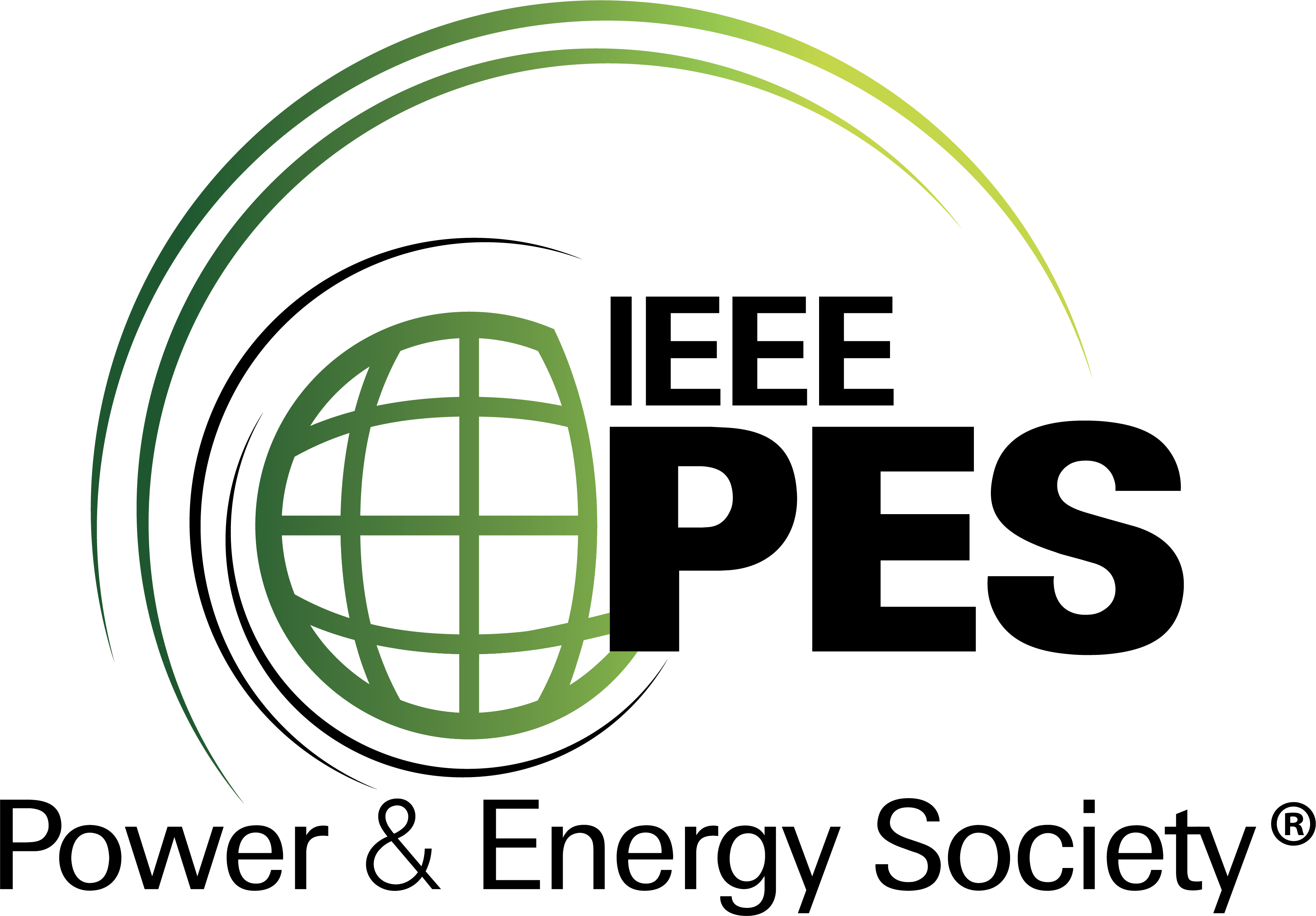}}

\usepackage{array}
\usepackage{stfloats}
\usepackage{url}
\usepackage{verbatim}
\usepackage{multirow}
\usepackage{xcolor}
\usepackage{gensymb}
\usepackage{algorithm}
\let\labelindent\relax % fix conflict between IEEEoj and enumitem
\usepackage{enumitem}
\usepackage{makecell}
\usepackage{xspace}
\usepackage{bm}
\usepackage{pgfplots}
\pgfplotsset{compat=1.17}

\newcommand{\bn}{retrospective }
\newcommand{\val}{CR-SE\xspace}
\newcommand{\pot}{Oracle}

\allowdisplaybreaks

\begin{document}

%\receiveddate{XX Month, XXXX}
%\reviseddate{XX Month, XXXX}
%\accepteddate{XX Month, XXXX}
%\publisheddate{XX Month, XXXX}
%\currentdate{XX Month, XXXX}
%\doiinfo{OAJPE.2023.1234567}

\title{Contextual Robust State Estimation in Distribution Systems {}{with Real-Time Unobservability and Scarce Data}}

\author{J. G. DE LA VARGA$^{1}$,
J. M. MORALES$^{1}$,
AND S. PINEDA$^{1}$,
}
\affil{OASYS Research Group, University of M\'alaga, Spain}
\corresp{CORRESPONDING AUTHOR: Juan M. Morales (e-mail: juan.morales@uma.es).}
\authornote{This work was supported by the Spanish Ministry of Science, Innovation and Universities (AEI/10.13039/501100011033) through project PID2023-148291NB-I00 and by the Department of Universities, Research and Innovation of the Regional Government of Andalusia through FEDER funds (grant PPRO-TEP967-G-2023). The work of J.~G.~De la Varga was supported by the Spanish Ministry of Science, Innovation and Universities training program for PhDs with fellowship number PRE2021-098958.}
\markboth{Contextual Robust State Estimation in Distribution Systems with Real-Time Unobservability and Scarce Data}{De la Varga \textit{et al.}}

\begin{abstract}
Distribution system state estimation faces a double information shortage: real-time measurements are too sparse to guarantee observability, \textcolor{black}{while historical data, however abundant in volume, is limited in its representativeness by the rapid evolution of network conditions}. This paper proposes a data-driven methodology that addresses both challenges simultaneously. We first introduce a $K$-nearest neighbor estimator that approximates the conditional distribution of delayed measurements given the available real-time context, formulating state estimation as a weighted least squares problem over plausible pseudo-measurement scenarios. Building on this, we develop a \emph{Contextual Robust State Estimator} (\val) that hedges against the inherent statistical uncertainty of the pseudo-measurements by optimizing over an adversarial reweighting of the nearest neighbors, subject to contextual proximity and monotonicity constraints. \textcolor{black}{\val hedges against the statistical uncertainty inherent to the generation of the pseudo-measurements---i.e., the potential misspecification of their conditional distribution---rather than against gross errors or bad data.} The resulting min-max problem admits an equivalent single-level reformulation of comparable tractability to the baseline estimator, and a data-driven procedure selects the robustness parameter online using the most recent training instance as a validation proxy. Numerical experiments on the IEEE 38-bus and a modified IEEE 123-bus \textcolor{black}{radial} networks\textcolor{black}{, as well as on the meshed IEEE 30-bus active distribution network,} demonstrate that \val strictly and consistently outperforms the baseline across all tested configurations of measurement availability, training size, and load level, with tail error improvements of up to 9\% and 11\% \textcolor{black}{in the radial networks and up to 18\% in the meshed one}, and the performance gap widening precisely in the configurations where the baseline is most vulnerable.
\end{abstract}


\maketitle

\section{INTRODUCTION}\label{sec:intro}

% {\color{orange}
% \IEEEPARstart{P}{ower} systems are transitioning toward decentralized architectures,
% driven by growing penetration of Distributed Energy Resources (DERs)
% whose stochastic generation and consumption introduce significant
% operational uncertainty. Active Distribution Networks (ADNs), with
% bidirectional flows and changing topology, make distribution-level
% operation and monitoring substantially more challenging than in the
% past~\cite{Aliwy2024ADN}.
% }

\IEEEPARstart{P}{ower} systems are undergoing a profound transformation, shifting from traditionally centralized architectures toward increasingly complex and decentralized networks. \textcolor{black}{The shift permeates every layer of system operation: generation, through integration with storage and flexible demand at the grid edge \cite{shaukat23}; voltage regulation, by migrating from centralized dispatch towards distributed and decentralized control schemes \cite{antoniadou17}; and monitoring, by increasing the need for surveillance over networks of cooperating agents, which raises in turn the question of how accurate state estimates can be sustained when some of those agents are compromised \cite{nasiri26}. These changes are} driven by the growing penetration of Distributed Energy Resources (DERs), such as rooftop solar, energy storage, and electric vehicles, which introduce significant operational uncertainty due to their stochastic behavior both in the generation and in the consumption. At the distribution level in particular, Active Distribution Networks (ADNs) characterized by bidirectional power flows, tighter interconnections and changing topology, emerge as a new paradigm, making their operation and monitoring substantially more challenging than in the past \cite{Aliwy2024ADN}.

Accurate state estimation in such uncertain environments remains a critical task for real-time monitoring and control. However, measurement infrastructure in distribution systems is typically sparse and heterogeneous \cite{yue2024graph}: fast but sparse measurements like SCADA systems and $\mu$-PMUs are received in real-time but exist in limited number, while slow but abundant measurements consist of more widely available data from Advanced Metering Infrastructure (AMI) such as Smart Meters (SM). This limitation in granularity and latency leads most existing approaches to favor static estimation frameworks over dynamic and forecasting-aided ones \cite{cheng2023survey}.

In the face of limited real-time measurements, historical data becomes a valuable resource for inferring the current state of the system. However, in modern distribution networks, even assuming that the grid topology and parameters are known, we typically face a \textit{double information shortage}. First, the available real-time measurements are often insufficient to guarantee system observability, especially given the sparse and uneven deployment of fast sensing infrastructure. Second, the historical data itself tends to be limited in both relevance and size, due to the rapid evolution of distribution systems: network topologies are frequently reconfigured and new DERs, such as rooftop solar units, are continuously being installed. As a result, previously collected data may quickly become outdated or uninformative. \textcolor{black}{Throughout this paper, \emph{scarce data samples} refers to the short window of recent instances collected under approximately the same operating regime---topology, DER configuration, and demand pattern---as the current one. The limitation is one of \emph{relevance} rather than volume: abundant but stale archives, such as billing or AMI histories, do not alleviate it.}

These two features, namely, real-time unobservability and \textcolor{black}{scarce data samples}, characterize the operational setting considered in this work. While several approaches to distribution system state estimation have addressed information limitations, a comprehensive treatment of the combined scenario involving real-time unobservability and \textcolor{black}{scarce data samples} remains largely unexplored in the literature, despite its growing practical relevance. In this paper, we propose estimation methods that are explicitly designed to operate under such information-constrained conditions, offering improved robustness and accuracy where conventional methods are not implementable or tend to fail.

Such a dramatic information shortage motivates the use of \emph{contextual learning} approaches: data-driven methods that use a \emph{context} vector (here, the real-time available measurements), which aim to exploit similarities between present and past measurements to approximate the conditional distribution of the underlying system state. However, this approximation is highly susceptible to inference errors, particularly when the available \textcolor{black}{data samples} are sparse or when the same real-time measurements may correspond to a broad range of possible system states~\cite{kumagai_constrained_2025}.

A prevalent strategy in the literature to handle previous challenges is the adoption of Bayesian frameworks. For instance, \cite{mestav_bayesian_2019} employs a deep Neural Network (NN) to learn mappings from limited real-time measurements to system states, combining it with a Bayesian state estimation algorithm to improve robustness against bad data. In \cite{dobbe_linear_2020}, Gaussian Process (GP)-based load forecasts are integrated with a linearized three-phase power flow model, enabling voltage predictions that are subsequently updated via linear Bayesian estimation. The authors in \cite{zhang_robust_2024} introduce a scalable and robust GP regression technique that supports accurate pseudo-measurement inference and state estimation using a limited number of measurements.

{ The authors of \cite{Donti2020} propose a matrix-completion--based state estimator tailored for low-observability distribution grids. Their method relies on a linearized power-flow model and is static, in the sense that it uses measurements from a single time instant. The approach exploits the empirically observed low-rank structure of the matrix collecting system variables to reconstruct missing information.

In a different line of research, several works incorporate data-driven corrections into physical models. For instance, \cite{liu_data_2025} introduces a data-enhanced linear power-flow formulation in which selected entries of the admittance matrix are adjusted using available measurements. Robustness is achieved through a Huber estimator, enabling convergence under model-parameter uncertainty, while zero-injection nodes are leveraged to improve observability. Similarly, \cite{zeraati_novel_2024} casts the problem as a mixed-integer linear program (MILP) in which both system states and nodal consumptions are decision variables. Pseudo-measurements derived from historical data define feasible bounds for load variables, and the MILP seeks to maximize state-estimation accuracy subject to power-flow constraints.} With a similar mathematical approach, \cite{karimi_joint_2021} addresses the challenge of topology identification in unobservable distribution systems with limited measurements. A compressive sensing framework that jointly estimates system state and network topology is proposed, formulated as a MILP with auxiliary variables to eliminate non-linearity and a convex relaxation to enable faster solution times. Finally, in \cite{pau_wls-based_2024} a modified weighted least squares (WLS) estimator is designed to operate effectively with very few real measurements and without relying on pseudo-measurements, even in typically unobservable scenarios. The method retains the advantages of a mature and widely used state estimation approach while enhancing observability by incorporating allocation factors, whose key assumption is that the power consumed (or generated) within each cluster of buses would be always at the same per-unit level (with respect to the rated one) for each node of the grid. \textcolor{black}{Neither the matrix completion of \cite{Donti2020} nor the modified WLS of \cite{pau_wls-based_2024} provides a conditional model of the delayed measurements given the real-time context, so neither is directly commensurable with the contextual setting considered here; the GP- and NN-based pseudo-measurement generators benchmarked in Section~\ref{sec:case_study}, in contrast, are.}

Despite their contributions, the described literature either assumes access to large amounts of historical data \cite{mestav_bayesian_2019,dobbe_linear_2020,zhang_robust_2024,liu_data_2025}, relies on approximated models \cite{Donti2020, liu_data_2025,zeraati_novel_2024,karimi_joint_2021}, or introduces structural simplifications that reduce the complexity of the state estimation \cite{Donti2020, pau_wls-based_2024}. Against this background, the contribution of this paper is twofold:

\begin{itemize}
\item Our work fits within the methodological line in which data-driven state estimation for non-observable distribution systems relies on the use of pseudo-measurements derived from \textcolor{black}{scarce data samples}. However, conventional approaches routinely incorporate these pseudo-measurements into weighted least-squares (WLS) estimators as if they were exact point predictions, even though they are statistical estimates obtained from finite samples. As a consequence, both the pseudo-measurements and their associated weights carry non-negligible statistical uncertainty that is systematically overlooked in the existing literature. \textcolor{black}{In contrast, this paper is, to the best of our knowledge, the first to explicitly robustify distribution system state estimation against the statistical misspecification of the conditional distribution of the pseudo-measurements. Moreover, the proposed robustification retains the computational complexity of a standard WLS estimator with pseudo-measurements, and the choice of its K-NN generator is supported empirically against GP- and NN-based alternatives rather than adopted for convenience.}

In this regard, it is worth contextualizing our contribution relative to the works \cite{zhang2023robust, vanin2023exact}, which are related but fundamentally distinct. In \cite{zhang2023robust}, the standard WLS estimator is extended to handle pseudo-measurements given as prediction intervals rather than point forecasts, with the goal of producing interval-valued state estimates. In \cite{vanin2023exact}, WLS is generalized to a maximum-likelihood framework to accommodate pseudo-measurements that follow non-Gaussian distributions. However, neither of these works addresses the crucial fact that both the prediction intervals in \cite{zhang2023robust} and the distributional assumptions in \cite{vanin2023exact} stem from data-driven models and therefore carry their own statistical uncertainty. Our method directly tackles this overlooked issue by providing robustness against misspecification of the pseudo-measurement distribution itself. \textcolor{black}{This notion of robustness is distributional and should not be confused with robustness against gross errors or adversarial interference on the real-time measurements. Therefore, our proposal is fully composable with any pre-estimation bad-data or false-data-injection detector.}

\item To construct the proposed robust estimator, we use a simple yet interpretable $K$-nearest neighbors (K-NN) mechanism to populate plausible pseudo-measurements as historical operating points whose real-time available measurements are closest to those observed at the estimation time (referred to as the \emph{context}). The proposed robust state estimator provides protection against pseudo-measurements' distributional misspecification without increasing the computational complexity of the underlying estimation problem. The resulting optimization model retains the same structure and tractability as a standard WLS estimator with pseudo-measurements, enabling straightforward implementation and real-time applicability.
% \textcolor{red}{Nevertheless, the KNN pseudo-measurement generator can be directly replaced by more advanced or application-specific pseudomeasurement-generation models---such as neural networks, linear or polynomial regressors, probabilistic forecasting tools, or any other pseudo-measurement generator---\emph{without modifying the core  robust formulation}}.
\end{itemize}

The rest of the paper is organized as follows. In Section~\ref{sec:methodology} we present the state estimation problem, emphasizing the challenges posed by unobservable systems and detailing the proposed approximation of the conditional distribution of delayed measurements given the real-time available ones. Both the non-robust baseline and the robust estimator are then presented, highlighting their respective advantages and limitations. Section~\ref{sec:case_study} provides numerical results that assess the performance of the proposed method against relevant benchmarks. Finally, Section~\ref{sec:conclusion} summarizes the key findings and outlines directions for future research.

\section{METHODOLOGY} \label{sec:methodology}

This section begins by presenting the conventional static state estimation problem in power systems, outlining its key components and highlighting a critical limitation of practical relevance: the lack of real-time observability due to the sparsity of available measurements. Building on this, we introduce two data-driven methods designed to estimate the system state under real-time unobservability and \textcolor{black}{scarce data samples}. Both methods rely on a $K$-nearest neighbors (K-NN) approach to approximate the conditional expectation of the residuals associated with delayed measurements, given the available real-time data. The second method (our proposal) incorporates a robustification mechanism aimed at mitigating inference errors caused by sample scarcity and the limited informativeness of the measurements. We then formalize this robust estimator mathematically and devote a separate section to detailing the data-driven procedure used to tune the robustness parameter that governs its behavior. 
% \textcolor{red}{ As discussed later, the use of K-NN is without loss of generality: it is adopted here due to its simplicity, reproducibility, and competitive performance, but it can be replaced by any other approach for conditional inference without altering the core methodology.}

Using generic notation, the relationship between the measurement vector $z \in \mathbb{R}^m$ and the state of the system $x \in \mathbb{R}^{n}$ is characterized as follows.
\begin{equation}\label{eq:z_x_relation}
{z}=h(x)+e
\end{equation}
where $h(\cdot):\mathbb{R}^{n} \rightarrow \mathbb{R}^{m}$ is the measurement function, which is constructed based on the power flow equations, the physical laws that govern the relationship between the measured magnitudes and the state variables. The term $e\in \mathbb{R}^{m}$ denotes the measurement error vector, whose components are typically assumed to follow normal distributions with zero mean and to be independent and uncorrelated with respect to the state variables \cite{abur2004}. Given the previous assumptions, the state estimation problem can be formulated as a WLS optimization problem \cite{abur2004}. This problem seeks to minimize the weighted sum of squares of the residuals:
\begin{subequations}
\label{eq:wls-se}
\begin{align}
\hat {x} \in \arg \min_x &\enskip (z - h(x))^{\top}W(z - h(x)) \\
\text{s.t.} & \enskip h^{0}(x)=0 \label{eq:wls_zero_inj}
\end{align}
\end{subequations}
where $\hat {x}$ is the estimated state vector that can be decomposed into voltage magnitude $\widehat{V}$ and angle $\hat{\theta}$, and the weight matrix $W=\text{diag}(\sigma^{-2}_1,\ldots,\sigma^{-2}_m)$ represents the user's confidence in the measured data. Here, $\sigma^2_m$ represents the error variance of the $m$-th measurement. Equation~\eqref{eq:wls_zero_inj} acts as a strong physical constraint at zero-injection nodes,
condition that enhances the observability of distribution systems and must be fully leveraged \cite{liu_data_2025}. For ease of notation, the solution of problem \eqref{eq:wls-se} is denoted as $\hat{x}=\texttt{wls}(z)$. In the general case, problem \eqref{eq:wls-se} results in a non-convex, constrained optimization problem, which has conventionally been solved using an iterative Newton-Raphson algorithm \cite{wang2019}. Alternatively, problem \eqref{eq:wls-se} can also be directly solved using interior-point methods implemented in non-linear optimization solvers~\cite{caro2020state}.

As discussed in Section \ref{sec:intro}, even when the number and types of measurements captured by the sensors are sufficient to make the system fully observable, some measurements may not be available in real time due to delays or communication issues.
In both this study and in prior work \cite{delaVarga2025}, we consider the measurement set divided into two groups: fast and sparse and slow but abundant. Consequently, we denote the set of available real-time measurements as $z^a\in\mathbb{R}^{m_a}$, the set of delayed measurements as $z^d\in\mathbb{R}^{m_d}$ { and, as shown in \eqref{eq:wls_zero_inj}, the set of zero-injection measurements $z^0\in\mathbb{R}^{m_0}$}. We then split the measurement function as follows:
\begin{subequations}\label{eq:za_zd_x_relation}
\begin{align}
& z^a=h^a(x)+e^a \\
& z^d=h^d(x)+e^d \\
& {z^0=0=h^0(x)}
\end{align}
\end{subequations}
where $h^a(\cdot):\mathbb{R}^{n} \rightarrow \mathbb{R}^{m_a}$, $h^d(\cdot):\mathbb{R}^{n} \rightarrow \mathbb{R}^{m_d}$, {$h^0(\cdot):\mathbb{R}^{n} \rightarrow \mathbb{R}^{m_0}$}, $e^a\in \mathbb{R}^{m_a}$, $e^d\in \mathbb{R}^{m_d}$, and $m=m_a+m_d{+m_0}$. In this work, we assume that the number and type of measurements available in real time are \emph{not} sufficient to make the system observable \cite{gomez-exposito2015, cheng2023survey}. 
% According to \cite{Baldwin1993ObservabilityPMU}, the system is numerically observable if the measurement Jacobian matrix is full rank and well-conditioned.

Although model \eqref{eq:wls-se} cannot be used to infer the state of the system in real time in the present instance $t$, denoted as $\hat{x}_t$, due to the lack of observability, we can always perform a \emph{\bn} state estimation for any past instance $s$ with measurements $(z^a_s,z^d_s)$. This retrospective system state is denoted as $\hat{x}^d_s$ and can be computed as $\hat{x}^d_s\in\texttt{wls}(z^a_s,z^d_s)$.  We assume that, when performing state estimation for the current instance $t$, we have access to the enlarged dataset $\{(z^a_s, z^d_s, \hat{x}^d_s),\  \forall s \in \mathcal{S}_t\}$, where $\mathcal{S}_t$ denotes the subset of $|\mathcal{S}_t|$ ordered \textcolor{black}{data samples} previous to present time $t$. In line with the online nature of the setting, we further assume that $\mathcal{S}_t$ has small cardinality, consistent with the idea that only the most recent past instances are informative for estimating the system state at time $t$. This choice directly addresses the tension between exploiting historical data and its limited representativeness in evolving systems: by restricting $\mathcal{S}_t$ to a short recent window, we confine the K-NN similarity argument to operating points collected under approximately the same topology, DER configuration, and demand regime as the current instance. \textcolor{black}{This short recent window is precisely the operational meaning of the \emph{scarce data samples} referred to in the title: not a shortage of archived volume, but a shortage of instances representative of current conditions.}

Although in practice we only have access to the available measurements \( z^a_t \), we can still estimate the conditional distribution of the delayed measurements \( z^d_t \) given \( z^a_t \) using a data-driven approach. Specifically, the
non-parametric K-NN method selects the \( K \) \textcolor{black}{data samples} whose available measurement vectors are closest to \( z^a_t \), identifying the corresponding set of indices \( \mathcal{S}_t^K \) as follows:
\begin{equation}
    \label{eq:K-NN}
    \mathcal{S}_t^K := \arg\min_{\substack{\mathcal{K} \subseteq \mathcal{S}_t \\ |\mathcal{K}| = K}} \sum_{k \in \mathcal{K}} {(z_t^a - z_k^a)^{\top} \hat{W}^a(z_t^a - z_k^a)}
\end{equation}
where $\hat{W}^a = \text{diag}(\hat{\sigma}^{-2}_1, \ldots, \hat{\sigma}^{-2}_{m_a})$ is a standardization matrix where $\hat{\sigma}_{m_a}$ denotes the empirical standard deviation of the $m_a$-th available measurement over the training set $\mathcal{S}_t$, ensuring that the distance metric is scale-invariant across measurements of heterogeneous nature and magnitude. \textcolor{black}{Note that $\hat{W}^a$ performs a per-measurement standardization, so that voltage magnitudes and power flows contribute to the distance on a common, dimensionless scale.} Consequently, the conditional distribution of the delayed measurements given the available ones, \( p(z^d_t \mid z^a_t) \), {is approximated by a uniform discrete distribution supported on the points \( \{ z^d_s : s \in \mathcal{S}_t^K \} \), each assigned a probability mass of \( 1/K \). Each of these points can be interpreted as a plausible realization of the missing measurements, i.e., a candidate pseudo-measurement vector.}

% \textcolor{red}{Importantly, as noted above, the K-NN mechanism used here to generate such scenarios is merely one possible choice: it can be directly replaced by any other conditional inference method capable of providing an estimate of the distribution of pseudo-measurements given the available data---such as neural networks, multivariate linear regression, or other probabilistic forecasting models---without requiring any modification to the core formulation of the state estimator we propose.}

In what follows, { we build upon the previous discussion and} present two data-driven methods for estimating the system state under real-time unobservability and \textcolor{black}{scarce data samples}. The first is a standard K-NN approach, which leverages nearest neighbors to approximate the conditional expectation of the delayed measurements' residuals given the available measurements { (i.e., to generate plausible pseudo-measurements)}. Although effective in many cases, this method is sensitive to noise and sample sparsity, inherent challenges in power system state estimation with limited data. To address this, we propose an extension inspired by distributionally robust optimization (DRO) with contextual information that explicitly accounts for uncertainty in the conditional inference process. By optimizing against ambiguity on the conditional distribution within a neighborhood of the one estimated by the K-NN, this robust formulation reduces inference errors, particularly in scenarios where data are scarce or non-representative.

\subsection{VANILLA K-NEAREST NEIGHBOR STATE ESTIMATOR} \label{sec:vanilla}
The first of the data-driven methods selects the $K$ \textcolor{black}{data samples} most similar to the real-time available measurements using \eqref{eq:K-NN} and searches for a state that best conforms with both the available measurements $z^a_t$ and the conditional probability distribution of the delayed measurements given the available measurements \( p(z^d_t \mid z^a_t) \) as estimated by the nearest neighbors.

Alternatively to the common approach of constructing pseudo-measurements \cite{delaVarga2025}, we adopt a \emph{stochastic} approach using plausible realizations or scenarios of the delayed measurements, making their associated uncertainty explicit in the formulation of the state estimator. This explicit treatment of uncertainty will later play a central role in enabling the robustification of our approach.

We denote this method \emph{vanilla} to differentiate it from the robust extension presented afterwards. The vanilla K-NN state estimate $\hat{x}^{v}_t$ is obtained by solving the non-linear weighted least squares problem \eqref{eq:vanilla}, initialized with a flat voltage. Instead of using fixed device error variances for the delayed measurements, we use the inverse of the sample variance of the $K$-nearest measurements in the training set, which indicates the confidence of the K-NN to predict its value. Analogously to the standard WLS formulation, a high variability across the neighbors is downweighted accordingly, while a measurement that is consistently similar across neighbors is treated with higher confidence. 
\begin{subequations}
\label{eq:vanilla}
\begin{align}
    \hat{x}_t^{v} \in
    \arg\min_{x} \
    &(z_t^a - h^a(x))^{\top} W^a (z_t^a - h^a(x))  \nonumber \\
    + \frac{1}{K} &\sum_{s \in \mathcal{S}_t^K} (z_s^d - h^d(x))^{\top}\hat{W}^d(z_s^d - h^d(x))  \\
    \text{s.t.} & \enskip h^{0}(x)=0
\end{align}
\end{subequations}

where $\hat{W}^d\!=\!\operatorname{diag}(\hat{\sigma}^{-2}_1, \ldots, \hat{\sigma}^{-2}_{m_d}) \text{, } (\hat{\sigma}^{d})^2\!=\!\tfrac{1}{K-1}\sum\nolimits_{s \in \mathcal{S}^K_t}(z^d_{s} - \bar{z}^d)^2 \text{ and } \bar{z}^d\!=\!\tfrac{1}{K}\sum\nolimits_{s \in \mathcal{S}^K_t} z^d_{s}$. \textcolor{black}{The weight $\hat{W}^d$ reflects the confidence of the conditional inference. Since it is computed over noisy meter readings $z^d_s$, it contains the device variance in addition to the conditional spread of the delayed measurements, and device precision further enters the pipeline through the retrospective states used as training targets.} This method is a \emph{lazy learning} approach, which means that no offline model training is required. Instead, the optimization problem is performed online using a local selection from the \textcolor{black}{scarce data samples}. The method is essentially equivalent to replacing the delayed and unavailable measurements with pseudo-measurements obtained from a K-NN prediction model.
The procedure is denoted $\hat{x}^{v}_t \in \texttt{Vanilla}(\mathcal{S}^K_t, z^{a}_{t})$, emphasizing its dependence on the $K$-nearest neighbors and the available measurements. 
% This approach is inherently sensitive to measurement noise and sample uncertainty, especially in the few-training-sample setting considered. To address this, we introduce a robust counterpart grounded in DRO with contextual information.

% \begin{algorithm}[h]
% \caption{\emph{Vanilla}} \label{alg:vanilla}
% \textbf{Input:} number of neighbors, $K$; training set, $\{(z_{s}^a, z_{s}^d) \enskip \forall \, s \in \mathcal{S}_t\}$; available measurements for test instance $z_{t}^a$.
% \begin{enumerate}[label={\arabic*)}]
%     \item Compute $\mathcal{S}_t^K$ using \eqref{eq:K-NN}.
%     \item Solve \eqref{eq:vanilla}: $\hat{x}_{t} = \hat{x}^{v}_{t} \in \texttt{Vanilla}(\mathcal{S}^K_t, z^{a}_{t})$.
% \end{enumerate}
% \textbf{Output:} State Estimation $\hat{x}_{t}$.
% \end{algorithm}

\subsection{CONTEXTUAL ROBUST STATE ESTIMATOR} \label{sec:robust}

The vanilla K-NN algorithm approximates the conditional distribution of the delayed measurement vector \( z^d_t \) given the available measurement vector \( z^a_t \) (referred to as the \emph{context}) using a discrete uniform distribution supported on the \( K \) training samples closest to \( z^a_t \). The underlying intuition is that training samples in close proximity to \( z^a_t \) are more likely to { serve as good pseudo-measurements in lieu of the missing vector \( z^d_t \). However, these  pseudo-measurements carry significant statistical uncertainty, particularly when the training dataset from which they have been obtained is small or in power systems where a single context \( z^a_t \) can correspond to a wide variety of system states}. To enhance robustness against such inherent statistical uncertainty, we draw inspiration from~\cite{esteban2022distributionally} and propose a simple yet effective \emph{ Contextual Robust State Estimator}, based on the min-max formulation~\eqref{eq:robust_K-NNse_minmax}.

\begin{subequations}
\begin{align}
    \hat{x}^{r}_{t} \in
    \mathop {\mathrm {arg\,min}}_{x}  \ \mathop {\mathrm {max}}_{w_s \geq 0} \
    & \sum_{s \in \mathcal{S}_t} w_s (z^{a}_{t} - h^a(x)) ^{\top} W^a(z^{a}_{t} - h^a(x))  \nonumber \\
    + \sum_{s \in \mathcal{S}_t} &w_s (z^d_{s} - h^d(x))^{\top}\hat{W}^d(z^d_{s} - h^d(x)) \\
     \text{s.t.} \enskip & w_s \leq \frac{1}{K}, \enskip\,\forall \, s \in \mathcal{S}_t \enskip :\lambda_s \label{eq:robust_K-NNse_minmax_K}\\
     & \sum_{s \in \mathcal{S}_t} w_s = 1 \enskip : \mu_1 \label{eq:robust_K-NNse_minmax_prob}\\
      w_{s+1} - & w_s \leq 0, \enskip\,\forall \, s \in \mathcal{S}_t \setminus \{s_{|\mathcal{S}_t|}\} \enskip :\gamma_s \label{eq:robust_K-NNse_minmax_weights}\\
      \sum_{s \in \mathcal{S}_t} w_s (z^{a}_{t} &- z^a_{s})^{\top}W^a(z^{a}_{t} - z^a_{s})  \leq \rho \enskip :\mu_2 \label{eq:robust_K-NNse_minmax_distance}\\
     & h^{0}(x)=0
    \end{align}\label{eq:robust_K-NNse_minmax}\\
\end{subequations}

In~\eqref{eq:robust_K-NNse_minmax}, the inner maximization accounts for the statistical uncertainty inherent in the K-NN-populated pseudo-measurements by selecting the most adversarial reweighting of the neighbors, subject to the four conditions defined by constraints~\eqref{eq:robust_K-NNse_minmax_K}--\eqref{eq:robust_K-NNse_minmax_distance}. Specifically, constraint~\eqref{eq:robust_K-NNse_minmax_K} ensures that no individual neighbor is assigned a weight greater than \( \frac{1}{K} \), which is equivalent to requiring that at least \( K \) neighbors receive non-zero weights. Constraint~\eqref{eq:robust_K-NNse_minmax_prob} enforces that the weights form a valid probability distribution with finite discrete support.  Constraint~\eqref{eq:robust_K-NNse_minmax_weights} imposes a monotonically decreasing structure on the weights to preserve the distance-based ranking of the neighbors, thus avoiding implausible conditional distributions where all the probability mass could be placed on the furthest neighbors. Lastly, constraint~\eqref{eq:robust_K-NNse_minmax_distance} bounds the squared, probability-weighted distance between the neighbors and the context \( z_t^a \) by a user-defined hyperparameter \( \rho \). Hence, this parameter plays the role of a tolerance margin that limits the allowable contextual deviation between the current available measurement vector $z^a_t$ and the available measurements \textcolor{black}{from the data samples} $z^a_s$. This ensures that the constructed probability distribution remains contextually close to $z^a_t$, while reducing the impact of the potential misalignment between similar observations and dissimilar underlying states. Expressed in the terminology of DRO, \( \rho \) controls the size of the ambiguity set and thereby governs the degree of reliance on the conditional distribution of the missing measurements estimated by the K-NN. In the jargon of state estimation, this conditional distribution characterizes plausible pseudo-measurements.
Therefore, through the inner maximization in~\eqref{eq:robust_K-NNse_minmax}, the proposed state estimator becomes \emph{insensitive} (up to the user-specified robustness radius $\rho$) to inaccuracies in the pseudo-measurements' weights. Intuitively, the maximization asks: ``What if the true distribution of pseudo-measurements differs from the empirical $1/K$ weights in the most adversarial way allowed by $\rho$? Will our estimation still be reliable?''

Finally, the outer minimization searches for a state estimate $x$ that conforms with both the available measurements and the least favorable reweighting of the neighbors (i.e., the plausible  pseudo-measurements vectors) under the four restrictions mentioned above. In doing so, the system state estimation is made robust against the unavoidable inference error incurred by the vanilla K-NN.   

Unlike the general distributionally robust approach with contextual information introduced in~\cite{esteban2022distributionally}, model~\eqref{eq:robust_K-NNse_minmax} possesses a distinct feature that makes it particularly appealing for enhancing robustness in the non-convex optimization involved in state estimation: it can be efficiently reformulated as the equivalent single-level minimization problem~\eqref{eq:robust_K-NNse}, which remains \emph{as tractable as} the vanilla K-NN state estimator~\eqref{eq:vanilla}. This reformulation requires only the addition of $2|\mathcal{S}_t| + 1$ variables and $3|\mathcal{S}_t|$ constraints (recall that the size of $\mathcal{S}_t$ is small) and leverages the fact that the inner maximization in~\eqref{eq:robust_K-NNse_minmax} is a linear program. As such, the transformation is achieved by deriving the dual of the inner maximization problem (with dual variables corresponding to each set of constraints in~\eqref{eq:robust_K-NNse_minmax} indicated after a colon) and combining the resulting minimization with the original outer problem as follows:

{
\begin{subequations}\label{eq:robust_K-NNse}
\begin{align}
\hat{x}^{r}_t \in &
\mathop{\mathrm {arg\,min}}_{x,\, \lambda_s,\, \gamma_s,\, \mu_1,\, \mu_2} \
 \frac{1}{K} \sum_{s \in \mathcal{S}_t} \lambda_s + \mu_1 + \rho\,\mu_2  \\
\text{s.t.} \
& \quad \lambda_s + \mu_1 + \mu_2\, (z^{a}_t - z^a_{s})^{\top}W^a(z^{a}_t - z^a_{s}) + f_s(\gamma) \geq \nonumber \\
 & \quad (z^{a}_t - h^a(x))^{\top}W^a(z^{a}_t - h^a(x))  \nonumber \\
+ & \quad (z^d_{s} - h^d(x))^{\top}\hat{W}^d(z^d_{s} - h^d(x)), \quad \forall s \in \mathcal{S}_t \\
& \quad \mu_2 \geq 0 \\
& \quad \lambda_s \geq 0, \quad \forall s \in \mathcal{S}_t \\
& \quad \gamma_s \geq 0, \quad \forall s \in \mathcal{S}_t \setminus \{s_{|\mathcal{S}_t|}\} \\
& \quad h^{0}(x)=0
\end{align}
\end{subequations}
}
where
$$f_s(\gamma) =
\begin{cases}
    -\gamma_s       & \text{if } s = s_1 \\
    \gamma_{s-1} - \gamma_s & \text{if } s_1 < s < s_{|\mathcal{S}_t|} \\
    \gamma_{s-1}    & \text{if } s = s_{|\mathcal{S}_t|}
\end{cases} \\$$

For a given $\rho$, we solve~\eqref{eq:robust_K-NNse} in two steps. We
first fix $x = \hat{x}^{v}_t$ (the vanilla estimate), which renders
the problem linear in the dual variables $\lambda_s, \gamma_s, \mu_1,
\mu_2$; solving this subproblem yields initial dual estimates. We then
solve the full non-linear problem over $(x, \lambda_s, \gamma_s,
\mu_1, \mu_2)$ warm-started from these initial estimates. We denote the resulting contextual robust estimate $\hat{x}^{r}_t \in
\texttt{Robust}(\mathcal{S}_t, z^{a}_{t}, K, \hat{x}^v_{t}, \rho)$.

\subsection{SELECTION OF THE ROBUSTNESS PARAMETER $\rho$} \label{sec:sel_rho}

The regularization parameter $\rho$ in \eqref{eq:robust_K-NNse} captures how informative the available measurements are about the delayed ones by limiting the weighted distance between the current vector $z_t^a$ and the \textcolor{black}{data samples} $z_s^a$, ensuring that only {close enough} instances from the past contribute to the estimation. Selecting an appropriate value of $\rho$ is therefore critical to achieving a good trade-off between robustness to inference errors and estimation accuracy.

To define a sensible data-driven range for the robustness parameter $\rho$,
we introduce two extreme values based on the structure of the problem. Since
constraint~\eqref{eq:robust_K-NNse_minmax_K} enforces that at least $K$
neighbors must receive non-zero weights, the feasibility of the problem
requires that $\rho$ be no smaller than the average contextual distance of the
$K$-nearest neighbors. This value thus defines a natural lower bound
$\rho_{min}$. Conversely, a loose upper bound $\rho_{max}$ can be set as the
contextual distance of the furthest neighbor in the dataset. These bounds are
formalized in \eqref{eq:rho_min}--\eqref{eq:rho_max}. In order to obtain a
grid of values spanning several orders of magnitude, we adopt a static set of
increments $\Delta\rho_c \, \forall \, c \in \mathcal{C}$, defined in \eqref{eq:delta_rhos}.
Consequently, candidate values are obtained as \eqref{eq:rho_sel} and any
candidate exceeding $\rho_{\max}$ is discarded, ensuring that the search
remains within the feasible range defined
by~\eqref{eq:rho_min}--\eqref{eq:rho_max}.
 
\begin{subequations} \label{eq:rhos}
\begin{gather}
    \rho_c = \rho_{\min} + \Delta\rho_c \leq \rho_{max} \quad \forall \, c \in \mathcal{C} \label{eq:rho_sel}\\
    \rho_{\min} = \frac{1}{K} \sum_{s \in \mathcal{S}_t^K}
     (z^{a}_{t} - z^a_{s})^{\top}W^a(z^{a}_{t} - z^a_{s})  \label{eq:rho_min} \\
    \rho_{\max} = \max_{s \in \mathcal{S}_t}\
    (z_t^a - z_s^a)^{\top} W^a(z_t^a - z_s^a) \label{eq:rho_max} \\
    \{\Delta\rho_c\}_{c \in \mathcal{C}} = \{0\} \cup \{e \times 10^{-f} : e \in \mathcal{E},\, f \in \mathcal{F}\}, \, \mathcal{E}, \mathcal{F} \subset \mathbb{Z}\label{eq:delta_rhos}
\end{gather}
\end{subequations}

Once the search range is defined, we perform an online grid search over the set of candidate values $\{\rho_c\,\forall\, c \in \mathcal{C}\}$ to obtain the corresponding state estimates $\hat{x}^r_{t,c}$. To select the optimal $\rho$, we exploit the fact that the most recent instance in the training set $\mathcal{S}_t$, denoted $s_{|\mathcal{S}_t|} \text{ or simply }t-1$, has a known retrospective state $\hat{x}^d_{t-1}$ and serves as a natural validation instance. Specifically, for each candidate $\rho_c$, we solve problem~\eqref{eq:robust_K-NNse} using the reduced training set $\mathcal{S}_t \setminus \{t-1\}$ and the available measurements of the previous instance $z^a_{t-1}$, in order to compute the estimation error against its retrospective state:
\begin{gather}
    \hat{x}^{r}_{t-1,c} \in \texttt{Robust}\left(\mathcal{S}_t\setminus\{t-1\}, z_{t-1}^a, K, \hat{x}^v_{t-1}, \rho_{c}\right) \, \forall \, c \in \mathcal{C} \\
    E_{t-1,c} = \|\hat{x}^r_{t-1,c} - \hat{x}^d_{t-1}\|^2_2
\end{gather}
The index $c^*$ that minimizes this validation error is then selected
for the final estimation on the actual test instance~$t$:
\begin{align} \label{eq:min_mse}
c^* &= \arg\min_{c \in \mathcal{C}} E_{t-1,c} \\
\hat{x}_t = \hat{x}^{r}_{t, c^*}&\in \texttt{Robust}\left(\mathcal{S}_t, z_{t}^a, K, \hat{x}^v_{t}, \rho_{c^*}\right)
\end{align}
This complete procedure, denoted \val (from \emph{Contextual Robust State Estimator}), is presented in Algorithm~\ref{alg:validation}.

\begin{algorithm}[h!]
\caption{\emph{Contextual Robust State Estimator} (\val )}\label{alg:validation}
\textbf{Input:} number of neighbors, $K$; measurements and \bn state in training set, $\{(z_{s}^a, z_{s}^d, \hat{x}^d_{s}) \enskip \forall \, s \in \mathcal{S}_t\}$;
set of robust parameter values, $\{\rho_{c} \enskip \forall \, c \in \mathcal{C}\}$.

\begin{enumerate}[label={\arabic*)}]
    \item Solve \eqref{eq:vanilla}: $\hat{x}^{v}_{t-1} \in \texttt{Vanilla}(\mathcal{S}^K_{t-1}, z^{a}_{t-1})$.
    \item Solve \eqref{eq:robust_K-NNse} : \\
    $\hat{x}^{r}_{t-1,c} \in \texttt{Robust}\left(\mathcal{S}_t\setminus\{t-1\}, z_{t-1}^a, K, \hat{x}^v_{t-1}, \rho_{c}\right) \, \forall \, c \in \mathcal{C}$.

    % Error
    \item Compute $E_{t-1,c} =
    ||\hat{x}^{r}_{t-1,c} - \hat{x}^d_{t-1}||_2
    \, \forall \, c \in \mathcal{C}$.

    \item Compute $c^* = \arg\min\limits_{c\in\mathcal{C}} E_{t-1,c}$.

    \item Solve \eqref{eq:vanilla}: $\hat{x}^{v}_{t} \in \texttt{Vanilla}(\mathcal{S}^K_t, z^{a}_{t})$.

    \item Solve \eqref{eq:robust_K-NNse} :
    $\hat{x}_{t} = \hat{x}^{r}_{t,c^*} \in \texttt{Robust}\left(\mathcal{S}_t, z_{t}^a, K, \hat{x}^v_{t}, \rho_{c^*}\right)$.
\end{enumerate}
\textbf{Output:} State Estimation $\hat{x}_{t}$
\end{algorithm}

\section{CASE STUDY} \label{sec:case_study}

In the following, we present results from a series of numerical experiments designed to evaluate the performance of our proposed methodology under varying operating conditions. We begin by outlining the experimental setup, which defines the test cases implemented, the time-series load profiles used and the selected low observability measurement configurations. This is followed by breaking down each of the estimators proposed in Section~\ref{sec:methodology} and comparing them against relevant benchmarks, detailing the metrics used. In the final part of the section, we discuss the observed results in depth, highlighting how each factor impacts performance and illustrating the consistency and benefits of the robust formulation across all scenarios.

\subsection{SETUP} \label{sec:setup}

We evaluate the proposed methodology on \textcolor{black}{three} distribution networks: the \textcolor{black}{radial} IEEE 38-bus system with 32 load nodes, the modified
IEEE 123-bus \textcolor{black}{radial} system of~\cite{bobo2021case123} with 85 load nodes\textcolor{black}{, and on the IEEE 30-bus system, a meshed network with distributed generation widely used in state estimation studies of modern, active systems~\cite{case30asADN_1, case30asADN_2}.}
\textcolor{black}{All networks} are modeled in Pandapower~\cite{pandapower2018}, and power flows are simulated \textcolor{black}{for the radial networks} using minute-level residential profiles from the Pecan Street dataset~\cite{pecanstreet}\textcolor{black}{, which already include behind-the-meter generation (predominantly rooftop PV)}.  We select day 100, which exhibits high intraday variability, making it a demanding test case. The load profiles are scaled to match the nominal loading of each network and then multiplied by a load factor $\delta \in \{0.5, 0.8, 1.0, 1.2\}$ spanning standard to heavily
stressed regimes, as quantified by the voltage statistics in Table~\ref{tab:voltage_stats}. \textcolor{black}{For the meshed system, the operating points are generated from an optimal power flow on the same time-series demands, which converges reliably only at light loading ($\delta=0.1$).} Power factors are held fixed at their nominal values. \textcolor{black}{The complete case-study data can be downloaded from \cite{noauthor_groupoasysrobustcontextualse_2025}.}

From the simulated day, we draw 100 test instances uniformly over
24~hours. For each, the training set $\mathcal{S}_t$ consists of its
$|\mathcal{S}_t|$ most recent predecessors, reflecting the online
nature of the setting. We consider $|\mathcal{S}_t| \in \{10, 20,
50\}$ and follow the commonly used heuristic $K\approx\sqrt{|\mathcal{S}_t|}$, yielding
$K \in \{3, 5, 7\}$. 
\textcolor{black}{A sensitivity sweep decoupling $K$ from $|\mathcal{S}_t|$ confirms the heuristic sits in a flat region of the estimation error: $K \in \{3,5,7\}$ perform within $\sim$1\% of each other, while $K=2$ and $K\geq10$ are 2--3\% worse.}
Measurements are obtained by adding Gaussian noise with $\sigma = 0.01$ p.u. for active and reactive power, and $\sigma = 0.001$ p.u. for voltage magnitudes, following
\cite{massignan2022bayesian}.

To study the effect of real-time measurement availability, we define
three configurations with an increasing number of SCADA measurements
$z^a$, detailed in Table~\ref{tab:redundancy_levels}. The most
constrained scenario includes only substation voltage magnitude and power flows; successive scenarios add line flow measurements (from-side) at selected branches, while preserving very low-observability conditions in real time. Smart Meter power
injections at load nodes are treated as delayed measurements $z^d$
due to slower reporting times~\cite{cheng2023survey}. Both $z^a$ and
$z^d$ are assumed fully accessible retrospectively, enabling training
and validation. 

\begin{table}
    \centering
    \caption{Mean and minimum nodal voltage magnitude (p.u.) across all buses and time steps, by network and load factor.}
    \label{tab:voltage_stats}
    \begin{tabular}{ccccc}
    \hline
    \multirow{2}{*}{$\delta$} & \multicolumn{2}{c}{38-bus} & \multicolumn{2}{c}{123-bus} \\
     & Mean & Min & Mean & Min \\
    \hline
    0.5 & 0.97 & 0.86 & 0.96 & 0.89 \\
    0.8 & 0.95 & 0.75 & 0.94 & 0.82 \\
    1.0 & 0.93 & 0.65 & 0.92 & 0.75 \\
    1.2 & 0.92 & 0.46 & 0.90 & 0.67 \\
    \hline
    \end{tabular}
\end{table}

\begin{table}[h]
    \centering
    \caption{Available measurements $z^a_t$ setup configurations.}
    \begin{tabular}{ccc}
        \hline
        \multirow{2}{*}{\centering\makecell{SCADA}} & \multicolumn{2}{c}{Network} \\
         & 38-bus & 123-bus \\
        \hline
        1 & $P, Q, V \,(\text{Substation})$  & $P, Q, V \,(\text{Substation})$  \\
        3 & $+ P_f, Q_f \, (\text{Lines }8, 12)$ & $+ P_f, Q_f \, (\text{Lines }6, 86)$ \\
        5 & $+ P_f, Q_f \, (\text{Lines }20, 32)$  & $+ P_f, Q_f \, (\text{Lines }9, 12)$  \\
        \hline
    \end{tabular}
    \label{tab:redundancy_levels}
\end{table}

\subsection{PERFORMANCE COMPARISON} \label{sec:comparison}
In this subsection, we introduce the learning-based estimators considered in our study and define the metrics used for performance comparison.
Thus, we define the following approaches:
\begin{itemize}

    \item [-] \emph{Vanilla}:
    This method corresponds to the vanilla K-NN state estimator described in Section~\ref{sec:methodology}-\ref{sec:vanilla}. It solves \eqref{eq:vanilla}, a WLS problem complemented with the $K$ closest delayed measurements \textcolor{black}{from the data samples} $z_s^d$ according to the available measurements $z_t^a$. No robustness or ambiguity modeling is applied in this method.

    \item [-] \val:
    This method is presented in Section~\ref{sec:methodology}-\ref{sec:robust} and corresponds to the contextual robust extension of \emph{Vanilla}. It leverages available measurements to estimate state variables in a manner that is robust to the uncertainty introduced by delayed measurements, which must be inferred from \textcolor{black}{the data samples}. This uncertainty arises from the variability in the relationship between available and delayed measurements, a challenge that becomes more pronounced when \textcolor{black}{data samples are} scarce.
    The full procedure is detailed in Algorithm~\ref{alg:validation}.

    \item [-] \emph{Retrospective}:
    This baseline assumes that all measurements are available and the system is then observable, that is, $\hat{x}_t = \hat{x}^d_t$. Under the considered setup, however, this situation is unrealistic and only used for benchmarking purposes.

    \item [-] \emph{\pot}:
    This baseline simulates an oracle version of \val,
    in which the \bn state $\hat{x}^d_t$ of the test instance is assumed to be known and available for hyperparameter tuning \textcolor{black}{(which only becomes available once the delayed measurements of instance $t$ arrive)}. While this benchmark is not implementable in practice, \textcolor{black}{it is reported solely to quantify the margin left by the data-driven selection of $\rho$}.

    \item [-] \emph{Persistent}:
    This baseline uses the state of the most recent training instance as the estimated state, that is, $\hat{x}_t = \hat{x}^d_{t-1}$. It avoids solving any online optimization problem, thus the computation is negligible. This method serves as a naive benchmark in the considered context of time-series forecasting, but is particularly relevant as this state is used in the validation process of the robust parameter $\rho$.
\end{itemize}

In order to compare the performance of the previous approaches, we use the Root Mean Squared Error (RMSE) defined as:
\begin{gather}
    \label{eq:rmse}
    \mathrm{RMSE} = \sqrt{\frac{1}{|\mathcal{T}|\,|\mathcal{N}|}\sum_{t \in \mathcal{T}}\sum_{i \in \mathcal{N}}E_{t,i}}, \quad E_{t,i} =|\hat{x}_{t,i} -\tilde{x}_{t,i}|^2
\end{gather}
where $\mathcal{T}$ and $\mathcal{N}$ are the sets of test instances and buses, $\hat{x}_{t,i}$ is the estimated state by each approach and $\tilde{x}_{t,i}$ is the ground-truth state, obtained from PF simulations as described in Section~\ref{sec:case_study}-\ref{sec:setup}.
In addition, as the distributional robustness of \val is explicitly designed to limit worst-case estimation errors, we also compute the Conditional Value at Risk (CVaR), which quantifies the average error in the tail of the distribution. As the considered distribution refers to the squared errors, analogously to \eqref{eq:rmse} we define the Root of the CVaR (RCVaR) as:
\begin{gather}
    \label{eq:cvar}
    \mathrm{CVaR}
    =
    \min_{\eta \in \mathbb{R}}
    \left\{
    \eta
    +
    \frac{1}{(1-\alpha)|\mathcal{T}|\,|\mathcal{N}|} \right. \nonumber \\
    \left. \sum_{t \in \mathcal{T}}\sum_{i \in \mathcal{N}}
    \max\left(
    E_{t,i} - \eta,\; 0
    \right)
    \right\} \\
    % E_t = \|\hat{x}_t - \tilde{x}_t\|_2^2 \\
    \mathrm{RCVaR} = \sqrt{\mathrm{CVaR}}
    \label{eq:rcvar}
\end{gather}
Finally, we select $\mathcal{E} = \{1, \ldots, 9\} \text{ and } \mathcal{F} = \{-1, \ldots, 2\}$, yielding $|\mathcal{C}|=37$ in \eqref{eq:rhos} and $\alpha=0.95$ in \eqref{eq:cvar}. State estimation problems \eqref{eq:vanilla} and \eqref{eq:robust_K-NNse} are modeled in Pyomo 6.8.0 \cite{pyomo} running in Python 3.11.10 and solved with Knitro 14.0.0 \cite{knitro}.

\subsection{RESULTS AND DISCUSSION}
To facilitate a detailed comparison of the methods introduced in Section~\ref{sec:case_study}-\ref{sec:comparison}, we begin by analyzing a representative scenario selected from the range of configurations discussed earlier. Specifically, Table~\ref{tab:results_base_case} reports estimation accuracy, measured as defined in~\eqref{eq:rmse} and \eqref{eq:cvar}, with 3 SCADA measurements (detailed in Table~\ref{tab:redundancy_levels}) and fixed conditions of moderate load level and training set size ($\delta = 0.5$, $|\mathcal{S}_t| = 20$, and $K = 5$), for both networks.

First, the results confirm that both the \emph{Retrospective} and the \emph{Persistent} baselines act accordingly as upper and lower bounds respectively. This highlights that relying exclusively on temporal correlations, that is, assuming the current state heavily mirrors the previous instance is insufficient to capture the real-time operational variations of distribution grids.

A direct comparison between the \emph{Vanilla} estimator and the proposed \val confirms the benefits of accounting for distributional uncertainty. In both distribution grids, \val strictly outperforms \emph{Vanilla}. Importantly, the \val estimator tracks the performance of the idealized \emph{\pot}, recovering 17-33\% in RMSE and 41-58\% in RCVaR of the gap between \emph{Vanilla} and the oracle method. This confirms that the validation procedure outlined in Section~\ref{sec:methodology}-\ref{sec:sel_rho} is effective for selecting the robustness parameter~$\rho$.

Table~\ref{tab:results_base_case} also illustrates that the gains from incorporating distributional robustness are similar for both systems in terms of mean errors, with relative improvements of around 2-2.5\% compared to \emph{Vanilla}. Notably, the introduction of robustness proves particularly effective at limiting extreme estimation error: in the smaller 38-bus network, the reduction in RCVaR is approximately 3.7\%, while in the larger 123-bus system, \val reduces the metric by approximately 6.7\%. This underscores that properly tuning the robustness parameter $\rho$ acts by constraining the influence of the most misleading pseudo-measurements, rather than uniformly shifting the estimated state.

\begin{table}
    \centering %

    \caption{Test set RMSE and RCVaR values ($\times 10^{-3}$) for all methods and both networks with 3 SCADA, $\delta=0.5,\,|\mathcal{S}_t|=20 \text{ and } K=5$.}
    \begin{tabular}{clcc}    % Tabla combinada
        \hline
       Network & Method & RMSE & RCVaR \\
         \hline
        \multirow{5}{*}{\centering\makecell{38-bus}}
        & \emph{Vanilla}	&	2.18	&	6.96	\\
        & \val	&	2.13	&	6.70	\\
        & \emph{Persistent}	&	4.45	&	14.62	\\
        & \emph{\pot}	&	2.03	&	6.51	\\
        & \emph{Retrospective} &	1.15	&	2.63	\\

         \hline
         \multirow{5}{*}{\centering\makecell{123-bus}}
        & \emph{Vanilla}	&	2.21	&	6.03	\\
        & \val	&	2.16	&	5.62	\\
        & \emph{Persistent}	&	3.69	&	9.85	\\
        & \emph{\pot}	&	1.91	&	5.03	\\
        & \emph{Retrospective} &	1.81	&	4.09	\\

         \hline
    \end{tabular}
    \label{tab:results_base_case}
\end{table}

\textcolor{black}{\textbf{Effect of the available measurements $z^a_t$.}} In what follows, we focus our analysis on comparing \emph{Vanilla} and \val methods.
The first block of Table~\ref{tab:results_sensitivity} evaluates the sensitivity of the estimators to the level of real-time observability by varying the number of available SCADA measurements from 1 to 5 (configured as in Table~\ref{tab:redundancy_levels}), maintaining moderate conditions of load level and training size ($\delta = 0.5$, $|\mathcal{S}_t| = 20$, and $K = 5$). The results of both the RMSE and RCVaR confirm that \val consistently outperforms the \emph{Vanilla} estimator across all evaluated measurement availability levels for both test networks. As expected, increasing SCADA availability reduces estimation errors for both methods across both networks. The two networks, however, exhibit contrasting sensitivity patterns that are informative about the nature of the ambiguity being addressed.

In the 38-bus network, the gains of \val over \emph{Vanilla} are largest under extreme unobservability: with only 1 SCADA measurement, \val outperforms \emph{Vanilla} by 5.9\% in RMSE and 6.7\% in RCVaR, a margin that narrows to 1.7\% in RMSE as observability improves to 5 SCADA measurements. This is consistent with the core premise of the approach: in a simple radial network, additional real-time measurements substantially reduce the ambiguity in the mapping between available and delayed variables, leaving less room for robustification to add value.

In the larger 123-bus network, the trend is reversed: the relative advantage of \val increases with additional measurements, with RMSE improvements rising from 2.3\% to 3.2\% and RCVaR improvements from 3.1\% to 6.0\% as SCADA availability grows from 1 to 5. This suggests that in a more complex system, limited local measurements cannot fully resolve system-wide ambiguity, making the optimal tuning of the robustness parameter $\rho$ essential for accurate global state estimation.

\begin{table}[t]
\setlength{\tabcolsep}{4pt}
\centering
\caption{\val vs.\ \emph{Vanilla (Vnl)} test-set RMSE and RCVaR ($\times10^{-3}$ p.u.) under three sensitivity sweeps. Each block varies a single axis with the others held fixed as indicated ($K\approx\sqrt{|\mathcal{S}_t|}$).}
\label{tab:results_sensitivity}
\begin{tabular}{l cccc cccc}
\hline
 & \multicolumn{4}{c}{38-bus} & \multicolumn{4}{c}{123-bus} \\
\cline{2-5}\cline{6-9}
 & \multicolumn{2}{c}{RMSE} & \multicolumn{2}{c}{RCVaR} & \multicolumn{2}{c}{RMSE} & \multicolumn{2}{c}{RCVaR} \\
 & \emph{Vnl} & \val & \emph{Vnl} & \val & \emph{Vnl} & \val & \emph{Vnl} & \val \\
\hline
\multicolumn{9}{l}{\emph{SCADA availability}\quad($\delta=0.5$, $|\mathcal{S}_t|=20$, $K=5$)} \\
\quad 1 & 3.62 & 3.41 & 11.72 & 10.94 & 2.38 & 2.32 & 6.57 & 6.37 \\
\quad 3 & 2.18 & 2.13 & 6.96 & 6.70 & 2.21 & 2.16 & 6.03 & 5.62 \\
\quad 5 & 2.13 & 2.10 & 6.73 & 6.61 & 2.15 & 2.08 & 5.88 & 5.52 \\
\hline
\multicolumn{9}{l}{\emph{Training size} $|\mathcal{S}_t|$\quad(1 SCADA, $\delta=0.5$)} \\
\quad 10 & 3.31 & 3.21 & 10.89 & 10.57 & 2.45 & 2.39 & 6.99 & 6.79 \\
\quad 20 & 3.62 & 3.41 & 11.72 & 10.94 & 2.38 & 2.32 & 6.57 & 6.37 \\
\quad 50 & 3.68 & 3.59 & 11.73 & 11.44 & 2.76 & 2.56 & 7.49 & 6.87 \\
\hline
\multicolumn{9}{l}{\emph{Load factor} $\delta$\quad(1 SCADA, $|\mathcal{S}_t|=50$, $K=7$)} \\
\quad 0.5 & 3.68 & 3.59 & 11.73 & 11.44 & 2.76 & 2.56 & 7.49 & 6.87 \\
\quad 0.8 & 6.27 & 5.87 & 20.85 & 19.07 & 4.30 & 3.95 & 12.50 & 11.13 \\
\quad 1.0 & 8.20 & 7.99 & 27.65 & 27.05 & 5.41 & 5.00 & 15.92 & 14.60 \\
\quad 1.2 & 10.05 & 9.64 & 33.51 & 32.53 & 6.45 & 6.12 & 19.19 & 17.96 \\
\hline
\end{tabular}
\end{table}

% \begin{table}
%     \caption{\val vs. \emph{Vanilla} test set RMSE and RCVaR values ($\times 10^{-3}$) for both networks with different SCADA availability, $\delta=0.5,\,|\mathcal{S}_t|=20 \text{ and } K=5$.}
%     \centering
%     \begin{tabular}{cccccc}  % SCADA + PMU
%         \hline
%         \multirow{2}{*}{Network}  & \multirow{2}{*}{SCADA} &  \multicolumn{2}{c}{RMSE} & \multicolumn{2}{c}{RCVaR} \\
%                             & & \emph{Vanilla} & \val & \emph{Vanilla} & \val \\
%         \hline
%         \multirow{3}{*}{\centering\makecell{38-bus}}
%         &	1	&	3.62	&	3.41	&	11.72	&	10.94	\\
%         &	3	&	2.18	&	2.13	&	6.96	&	6.70	\\
%         &	5	&	2.13	&	2.10	&	6.73	&	6.61	\\

%         \hline
%         \multirow{3}{*}{\centering\makecell{123-bus}}
%         &	1	&	2.38	&	2.32	&	6.57	&	6.37	\\
%         &	3	&	2.21	&	2.16	&	6.03	&	5.62	\\
%         &	5	&	2.15	&	2.08	&	5.88	&	5.52	\\

%         \hline
%     \end{tabular}
%     \label{tab:results_measurements}
% \end{table}

\textcolor{black}{\textbf{Effect of the training set $|\mathcal{S}_t|$.}} The influence of the training set size across network types is assessed in the second block of Table~\ref{tab:results_sensitivity}. We continue considering moderate loading  and we restrict the analysis to the 1 SCADA configuration, the most common and demanding at the distribution level. The degradation in the \emph{Vanilla} estimator as the \textcolor{black}{data sample} window grows confirms a known theoretical vulnerability of K-NN regression. With minimal $z^a$, expanding $|\mathcal{S}_t|$ from 10 to 50 worsens \emph{Vanilla}'s RMSE by roughly 11--12\% and the RCVaR by 7\%. This occurs because, with minimal real-time context, a larger pool of candidates increases the likelihood of retrieving historically distant operating points that happen to match current SCADA readings. This increases the uncertainty in the delayed pseudo-measurements, which in turn harms the state estimation.

The \val estimator substantially attenuates this degradation: by tuning $\rho$ to penalise contextually distant neighbors, it effectively shrinks the influence of ambiguous pseudo-measurements without discarding them entirely. Crucially, the performance gap between \val and \emph{Vanilla} widens as $|\mathcal{S}_t|$ increases. In the 123-bus network, the relative RCVaR improvement grows from 2.9\% ($|\mathcal{S}_t|=10$) to 8.2\% ($|\mathcal{S}_t|=50$), suggesting that the robustification mechanism scales beneficially with data abundance---a desirable property for practical deployment where training archives naturally grow over time.

% \begin{table}
% \centering
% \caption{\val vs. \emph{Vanilla} test set RMSE and RCVaR values ($\times 10^{-3}$) for both networks with 1 SCADA and $\delta=0.5$.}
% \label{tab:results_training_size}
% \begin{tabular}{ccccccc}
% \hline
% \multirow{2}{*}{Network} & \multirow{2}{*}{$|\mathcal{S}_t|$} & \multirow{2}{*}{$K$} & \multicolumn{2}{c}{RMSE} & \multicolumn{2}{c}{RCVaR} \\
% & & & \emph{Vanilla} & \val & \emph{Vanilla} & \val \\
% \hline
% \multirow{3}{*}{38-bus}

% &	10	&	3	&	3.31	&	3.21	&	10.89	&	10.57	\\
% &	20	&	5	&	3.62	&   3.41    &   11.72   &   10.94	\\
% &	50	&	7	&	3.68	&	3.59	&	11.73	&	11.44	\\

% \hline
% \multirow{3}{*}{123-bus}
% &	10	&	3	&	2.45	&	2.39	&	6.99	&	6.79	\\
% &	20	&	5	&	2.38	&	2.32	&	6.57	&	6.37	\\
% &	50	&	7	&	2.76	&	2.56	&	7.49	&	6.87	\\

% \hline
% \end{tabular}
% \end{table}

\textcolor{black}{\textbf{Effect of the load factor $\delta$.}} Lastly, the effect of load demand variation is studied in the third block of Table~\ref{tab:results_sensitivity}. As the load factor $\delta$ increases from 0.5 to 1.2, both RMSE and RCVaR grow substantially for all methods, with errors roughly tripling between the lowest and highest load levels. This is expected: heavier and more variable loading induces stronger nonlinearities in the power flow equations, making it harder to infer the system state from a sparse set of real-time measurements.

The relative performance of \val over \emph{Vanilla} does not follow a direct correlation with $\delta$. In the 38-bus network, the largest relative RCVaR gain occurs at $\delta= 0.8$ (8.5\%), while the improvement is more modest at the rest of load levels (2--3\%). In the 123-bus network, gains are more uniformly distributed across the first three load levels (8-11\%), before declining slightly at the heaviest loading (6.4\%). In any case, \val strictly outperforms \emph{Vanilla} across all tested load levels and both networks, confirming that the robustification provides consistent benefits under diverse operating conditions.

% \begin{table}
%     \caption{\val vs. \emph{Vanilla} test set RMSE and RCVaR values ($\times 10^{-3}$) for both networks and different load levels $\delta$, with 1 SCADA, $|\mathcal{S}_t|=50 \text{ and } K=7$.}
%     \centering
%     \begin{tabular}{cccccc}  % SCADA + PMU
%         \hline
%         \multirow{2}{*}{Network}  & \multirow{2}{*}{Load Factor} &  \multicolumn{2}{c}{RMSE} & \multicolumn{2}{c}{RCVaR} \\
%                             & & \emph{Vanilla} & \val & \emph{Vanilla} & \val \\
%         \hline
%         \multirow{4}{*}{\centering\makecell{38-bus}}
%         &	0.5	&	3.68	&	3.59	&	11.73	&	11.44	\\
%         &	0.8	&	6.27	&	5.87	&	20.85	&	19.07	\\
%         &	1	&	8.20	&	7.99	&	27.65	&	27.05	\\
%         &	1.2	&	10.05	&	9.64	&	33.51	&	32.53	\\

%         \hline
%         \multirow{4}{*}{\centering\makecell{123-bus}}
%         &	0.5	&	2.76	&	2.56	&	7.49	&	6.87	\\
%         &	0.8	&	4.30	&	3.95	&	12.50	&	11.13	\\
%         &	1	&	5.41	&	5.00	&	15.92	&	14.60	\\
%         &	1.2	&	6.45	&	6.12	&	19.19	&	17.96	\\

%         \hline
%     \end{tabular}
%     \label{tab:results_load}
% \end{table}

\textcolor{black}{\textbf{Statistical analysis.} Beyond averages, further information can be retrieved from the statistics of the improvement in Table~5. At the level of individual instances, \val is marginally worse than \emph{Vanilla} only on the most benign instances, while being significantly better in aggregate on both metrics: the robustification deliberately trades small losses on easy instances for large error reductions on the damaging ones. This is the intended behavior of a robust estimator, and it explains why the RCVaR improves consistently more than the RMSE---by a factor of 1.3 on the 38-bus network and 1.6 on the 123-bus network. Operationally, it is these tail errors and not the average, the ones that drive voltage-limit violations and mis-informed control actions, so a reduction concentrated in the tail is the most relevant kind for real-time operation. Consistently, the gains are largest where the baseline is most exposed: under extreme unobservability, the 38-bus RMSE improvement rises to $+3.7\%$ with 1 SCADA measurement, against $+0.8$--$1.0\%$ with 3--5.}

\begin{table}[t]
\setlength{\tabcolsep}{3pt}
{\color{black}
\centering
\caption{\val vs.\ \emph{Vanilla} over the full grid of 36 configurations per network: improvement (\%), with an exact binomial sign test across configurations.}
\label{tab:Rsig}
\begin{tabular}{llcccc}
\hline
Network & Metric & Mean & Median & Configs.\ improved & Sign-test $p$ \\
\hline
\multirow{2}{*}{38-bus}  & RMSE  & $+1.87$ & $+2.39$ & 27 / 36 & $2.0\times10^{-3}$ \\
                         & RCVaR & $+2.36$ & $+2.34$ & 30 / 36 & $3.5\times10^{-5}$ \\
\hline
\multirow{2}{*}{123-bus} & RMSE  & $+2.22$ & $+2.26$ & 26 / 36 & $5.7\times10^{-3}$ \\
                         & RCVaR & $+3.73$ & $+4.44$ & 29 / 36 & $1.6\times10^{-4}$ \\
\hline
\end{tabular}}
\end{table}

\begin{figure}[t]
    \centering
    \includegraphics[width=1\linewidth]{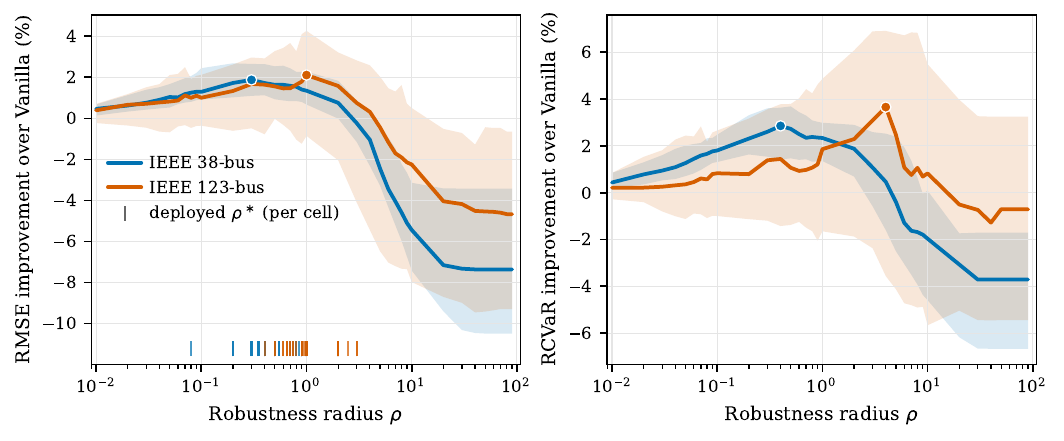}
    \caption{\textcolor{black}{Improvement of \val over \emph{Vanilla} when $\rho$ is held fixed across its candidate range, for both networks (left: RMSE; right: RCVaR). Solid line: median across the 36 configurations per network; shaded: inter-quartile range. Markers indicate the median peak. The rug at the bottom of the left panel shows the $\rho$ deployed by Algorithm~\ref{alg:validation} (median across test instances per configuration), which lands inside the beneficial plateau. $\rho = 0$ recovers \emph{Vanilla} and is omitted from the log axis.}}
    \label{fig:rho_sensitivity}
\end{figure}

\textcolor{black}{\textbf{Effect of the hyperparameter $\rho$.} Fig.~\ref{fig:rho_sensitivity} shows the improvement of the robust estimator over \emph{Vanilla} when $\rho$ is held fixed across its candidate range, for all 100 instances of the 36 configurations. The response is a broad plateau rather than a knife edge: the beneficial range spans three orders of magnitude, from the lower end of the candidate grid up to $\rho \approx$ 2--4. Overly large values of $\rho$ are uniformly harmful, showing the expected over-hedging behavior. Peak values on the 38-bus network appear at $+1.88\%$ at $\rho=0.3$ for the RMSE ($+2.86\%$ at $\rho=0.4$ for the RCVaR), while on the 123-bus network at $+2.11\%$ at $\rho=1$ for the RMSE ($+3.66\%$ at $\rho=4$ for the RCVaR). In 67 of the 72 configurations there exists some $\rho>0$ that improves on \emph{Vanilla}, and the values selected by the validation procedure of Section~\ref{sec:methodology}-\ref{sec:sel_rho} (whose median values appear on the rug at the bottom of the left panel) fall inside this plateau on both networks, which explains why the data-driven selection recovers a substantial share of the achievable gain (Table~\ref{tab:results_base_case}).}

\textcolor{black}{\textbf{Comparison with alternative pseudo-measurement generators.} To assess whether more expressive conditional models outperform the proposed estimator in the regime under study, we benchmark two alternative generators of the delayed pseudo-measurements: a Gaussian Process Regressor (GPR), fitted independently per delayed bus, and a Deep Neural Network (DNN), regressing all delayed buses jointly. Both are embedded in the \emph{same} WLS estimator, fitted on the \emph{same} training windows $\mathcal{S}_t$, and weighted the \emph{same} way, so that only the pseudo-measurement generator differs. Due to space limitations, their full configuration is provided with the case-study data~\cite{noauthor_groupoasysrobustcontextualse_2025}.}

\textcolor{black}{Table~\ref{tab:Rbench} reports the improvement of \val over each benchmark across the 36 configurations of each network. \val attains a lower RMSE than every benchmark on both networks. Against the strongest parametric contender, the GPR, the margin is 3.6\% (38-bus) and 4.9\% (123-bus) in RMSE and 4.3\% / 5.5\% in RCVaR, significant in both cases under an exact sign test across configurations. The DNN, however, being the most flexible model, is essentially overfitting training data and performs worst.}

\begin{table}[t]
{\color{black}
\setlength{\tabcolsep}{3pt}
\centering
\caption{Mean relative improvement of \val over each pseudo-measurement generator benchmark over the 36 configurations per network (positive $=$ \val better; in parentheses, configurations in which \val is better, out of 36). All comparisons significant under an exact binomial sign test ($p<0.05$).}
\label{tab:Rbench}
\small
\begin{tabular}{lcccc}
\hline
 & \multicolumn{2}{c}{RMSE (\%)} & \multicolumn{2}{c}{RCVaR (\%)} \\
\cline{2-3}\cline{4-5}
Benchmark & 38-bus & 123-bus & 38-bus & 123-bus \\
\hline
GPR         & $+3.58$ (31) & $+4.87$ (32) & $+4.27$ (29)  & $+5.50$ (27) \\
DNN                & $+9.03$ (33) & $+19.95$ (35) & $+9.28$ (28) & $+25.88$ (36) \\
K-NN (\emph{Vanilla})      & $+1.87$ (27) & $+2.22$ (26) & $+2.36$ (30)  & $+3.73$ (29) \\
\hline
\end{tabular}}
\end{table}

\textcolor{black}{\textbf{A meshed system: the IEEE 30-bus active distribution network.} 
As introduced in Section~\ref{sec:case_study}-\ref{sec:setup}, for this system we select a single load level $\delta=0.1$ because, as studied in the radial systems, the load factor is the experimental axis the improvement is least sensitive to. Across the 9 tested configurations of Table~\ref{tab:case30} (3 SCADA levels $\times$ 3 $|\mathcal{S}_t|$ levels), \val improves on \emph{Vanilla} in every one (sign test $p=0.002$), with a mean RMSE improvement of 7.0\% (median 3.6\%) and a mean RCVaR improvement of 7.4\% (median 3.2\%). The largest gains (among all three networks) are concentrated at the smallest training window (12.6--18.1\% in RMSE, 12.5--18.1\% in RCVaR), consistent with the scarce-data motivation of the method. We attribute these larger gains to the stronger contextual ambiguity of an active, meshed network: distributed generation and a meshed topology make a single real-time context consistent with a wider range of nodal states, and hence with more dispersed pseudo-measurements. Since the robust layer hedges precisely against this inference uncertainty, its benefit is largest exactly where that uncertainty is greatest.}

\begin{table}[t]
{\color{black}
\centering
\caption{\val on the IEEE 30-bus active distribution network, for $\delta=0.1$. RMSE and RCVaR in $10^{-3}$ p.u.}
\label{tab:case30}
\begin{tabular}{c|cc|cc|cc}
\hline
& \multicolumn{2}{c|}{1 SCADA} & \multicolumn{2}{c|}{3 SCADA} & \multicolumn{2}{c}{5 SCADA} \\
$|\mathcal{S}_t|$ & \emph{Vanilla} & \val & \emph{Vanilla} & \val & \emph{Vanilla} & \val \\
\hline
\multicolumn{7}{c}{RMSE} \\ \hline
10 & 3.83 & 3.35 & 3.99 & 3.27 & 3.63 & 3.12 \\
20 & 3.18 & 3.12 & 3.22 & 3.10 & 3.31 & 3.09 \\
50 & 3.38 & 3.30 & 3.39 & 3.31 & 3.31 & 3.24 \\
\hline
\multicolumn{7}{c}{RCVaR} \\ \hline
10 & 10.54 & 9.22 & 11.15 & 9.14 & 10.15 & 8.47 \\
20 & 8.86 & 8.78 & 9.08 & 8.79 & 9.32 & 8.60 \\
50 & 9.80 & 9.57 & 9.82 & 9.51 & 9.17 & 8.96 \\
\hline
\end{tabular}}
\end{table}

\subsection{PRACTICAL CONSIDERATIONS}
\label{sec:limitations}

\textcolor{black}{\textbf{Scalability.} The robust reformulation~\eqref{eq:robust_K-NNse} adds $2|\mathcal{S}_t|+1$ variables and $3|\mathcal{S}_t|$ constraints to the vanilla WLS problem---overhead that scales with the (small) training size and is independent of the network size. A single robust solve takes on average 3.7--4.3~s on the 123-bus network, against 3.8--4.0~s for a vanilla solve (on the 38-bus network, 0.6--0.8~s versus $\sim$0.5~s); in solver effort, the robust problem requires 1.4--2$\times$ the interior-point iterations of the vanilla one, independently of network size. The complete procedure of Algorithm~\ref{alg:validation} solves the robust problem once per candidate $\rho$ on the validation instance, but the candidate solves are independent and can be parallelized. Scaling to larger feeders is therefore governed by the base non-linear WLS solve, common to all WLS-based estimators, rather than by the robustification layer. Additionally, we report that the solve time is essentially flat in $K$: across the pairs $(|\mathcal{S}_t|,K)\in\{(10,3),(20,5),(50,7)\}$ it varies by less than 15\%, a growth attributable to the $O(|\mathcal{S}_t|)$ added dual variables and constraints rather than to $K$, which enters only through the bound $1/K$ in \eqref{eq:robust_K-NNse_minmax_K} and does not change the problem dimension.}

\textcolor{black}{\textbf{Topology variability.} The estimator assumes the network topology (and hence $h$) is fixed across $\mathcal{S}_t$ and the test instance, and two cases must be distinguished. 
When the switching status is directly monitored by SCADA breakers, a reconfiguration is observed and $\mathcal{S}_t$ can be conditioned on the current topology, excluding any instance collected under a different switch state.
When the reconfiguration is unmonitored, the change in $h$ is latent and only partially mitigated by the deliberately small window: at minute resolution, $|\mathcal{S}_t|$ spans from only ten minutes to at most about one hour.
Removing the assumption altogether, by jointly estimating $h$ through topology identification or an additional robustification layer over the measurement function, is a natural extension.}
% \textbf{Topology variability.} The estimator assumes the network
% topology (and hence $h$) is fixed across $\mathcal{S}_t$ and the test
% instance. Reconfiguration within this window would partially break the K-NN similarity argument, a risk mitigated by the small
% $|\mathcal{S}_t|$. Alternatively, instances collected under topologies dissimilar to the current one can be filtered out from $|\mathcal{S}_t|$ prior to estimation, effectively restricting the training set to topologically consistent operating points. Jointly estimating $h$---via topology identification or an additional robustification layer---is a natural extension.

\textcolor{black}{\textbf{Unbalanced operation.} The balanced single-phase modeling is a limitation of the present study. The K-NN generator and the robust layer act on the pseudo-measurement distribution, not on the network model, and therefore transfer to a three-phase formulation without reformulation; only the measurement function $h$ and the dimensions of $(x,z)$ change, at the cost of a larger base WLS solve and of a per-phase correlation structure in the delayed measurements. A full unbalanced study, e.g., via a three-phase power flow, is a concrete next step.}

% \textbf{Adversarial perturbations.} The ambiguity set hedges against
% \emph{statistical} misspecification of the pseudo-measurements, not
% adversarial interference. While DRO may incidentally dampen small
% biased perturbations, false data injection on the real-time channel
% itself requires dedicated bad-data detection or adversarial-robust
% formulations that are orthogonal to---and composable with---the one
% developed here.
\textbf{Bad data and adversarial perturbations.} The ambiguity set
hedges against \emph{statistical} misspecification of the
pseudo-measurements, not against gross errors or
adversarial interference on the real-time measurements. Classical post-estimation bad-data detection and identification (BDDI) based on residual analysis is computationally demanding, as it typically requires multiple re-runs of the state estimator. Pre-estimation BDDI of corrupted measurements, followed by estimation on the cleaned input, is therefore both cheaper and more reliable~\cite{Kfouri2026bddi}. Our formulation is fully composable with this preferred pipeline: any pre-estimation bad-data or false-data-injection detector can be applied upstream, with \val\ operating on the sanitized $(z^a_t, z^d_s)$. 
% Developing a detector tailored to the low-observability regime considered here is orthogonal to---and complementary with---the contribution developed in this paper.

\section{CONCLUSION}
\label{sec:conclusion}

This work addresses state estimation in distribution systems under a
\emph{double information shortage}: partial unobservability from
scarce real-time measurements, and \textcolor{black}{scarce data samples} from
evolving network infrastructure. The proposed \val combines a K-NN
pseudo-measurement generator with a contextual DRO robustification
layer to reduce the sensitivity to potentially misleading \textcolor{black}{data samples}, and admits an equivalent single-level reformulation of
comparable tractability to standard WLS with similar pseudo-measurements.

Numerical experiments on the IEEE 38-bus and IEEE 123-bus networks\textcolor{black}{, together with a meshed IEEE 30-bus active distribution network,}
confirm that \val strictly and consistently outperforms the
\emph{Vanilla} K-NN baseline across all tested configurations, with
relative improvements of up to 6.41\% in RMSE and 8.54\% in RCVaR on
the 38-bus system, and up to 8.09\% and 10.95\% respectively on the
123-bus system. \textcolor{black}{On the meshed IEEE 30-bus network, where contextual ambiguity is greatest, the gains are largest, reaching up to 18.1\% in both RMSE and RCVaR at the smallest training window.} Gains are most pronounced where \emph{Vanilla} is most
vulnerable: under extreme unobservability, with larger training sets,
and at heavier loading levels where nonlinearities widen the
conditional uncertainty of the delayed measurements. These results
validate both the formulation and the data-driven procedure for
selecting~$\rho$, which not only improves average accuracy but
specifically targets the most damaging estimation errors.

\section*{ACKNOWLEDGMENT}
The authors thankfully acknowledge the computer resources (Picasso Supercomputer), technical expertise and assistance provided by the SCBI (Supercomputing and Bioinformatics) center of the University of M\'alaga.

\bibliographystyle{IEEEtran}
\bibliography{references}

\begin{thebibliography}{10}
\providecommand{\url}[1]{#1}
\csname url@samestyle\endcsname
\providecommand{\newblock}{\relax}
\providecommand{\bibinfo}[2]{#2}
\providecommand{\BIBentrySTDinterwordspacing}{\spaceskip=0pt\relax}
\providecommand{\BIBentryALTinterwordstretchfactor}{4}
\providecommand{\BIBentryALTinterwordspacing}{\spaceskip=\fontdimen2\font plus
\BIBentryALTinterwordstretchfactor\fontdimen3\font minus \fontdimen4\font\relax}
\providecommand{\BIBforeignlanguage}[2]{{%
\expandafter\ifx\csname l@#1\endcsname\relax
\typeout{** WARNING: IEEEtran.bst: No hyphenation pattern has been}%
\typeout{** loaded for the language `#1'. Using the pattern for}%
\typeout{** the default language instead.}%
\else
\language=\csname l@#1\endcsname
\fi
#2}}
\providecommand{\BIBdecl}{\relax}
\BIBdecl

\bibitem{shaukat23}
N.~Shaukat, M.~R. Islam, M.~M. Rahman, B.~Khan, B.~Ullah, S.~M. Ali, and A.~Fekih, ``Decentralized, democratized, and decarbonized future electric power distribution grids: A survey on the paradigm shift from the conventional power system to micro grid structures,'' \emph{IEEE Access}, vol.~11, pp. 60\,957--60\,987, 2023.

\bibitem{antoniadou17}
K.~E. Antoniadou-Plytaria, I.~N. Kouveliotis-Lysikatos, P.~S. Georgilakis, and N.~D. Hatziargyriou, ``Distributed and decentralized voltage control of smart distribution networks: Models, methods, and future research,'' \emph{IEEE Transactions on Smart Grid}, vol.~8, no.~6, pp. 2999--3008, Nov. 2017.

\bibitem{nasiri26}
S.~Nasiri, H.~Seifi, and H.~Delkhosh, ``A trust-aware consensus mechanism for post-attack restoration of power system distributed state estimation,'' \emph{IEEE Transactions on Smart Grid}, vol.~17, no.~1, pp. 562--575, 2026.

\bibitem{Aliwy2024ADN}
A.~K. Aliwy, M.~Al-badri, N.~Q. Ali, H.~F.~K. Al-khazaali, and D.~M. Hussein, ``Current state and future trends of active management in distribution networks,'' in \emph{2024 21st International Multi-Conference on Systems, Signals \& Devices (SSD)}, 2024, pp. 766--775.

\bibitem{yue2024graph}
H.~Yue, W.~Zhang, U.~C. Yilmaz, T.~Yildiz, H.~Huang, H.~Liu, Y.~Lin, and A.~Abur, ``Graph-learning-assisted state estimation using sparse heterogeneous measurements,'' \emph{Electric power systems research}, vol. 235, p. 110644, 2024.

\bibitem{cheng2023survey}
G.~Cheng, Y.~Lin, A.~Abur, A.~Gómez-Expósito, and W.~Wu, ``A survey of power system state estimation using multiple data sources: Pmus, scada, ami, and beyond,'' \emph{IEEE Transactions on Smart Grid}, vol.~15, no.~1, pp. 1129--1151, 2024.

\bibitem{kumagai_constrained_2025}
M.~Kumagai and M.~Watanabe, ``Constrained state estimation with probabilistic error model for underdetermined systems such as distribution network,'' \emph{IEEJ Transactions on Electrical and Electronic Engineering}, vol.~20, no.~8, pp. 1186--1194, 2025.

\bibitem{mestav_bayesian_2019}
K.~R. Mestav, J.~Luengo-Rozas, and L.~Tong, ``Bayesian state estimation for unobservable distribution systems via deep learning,'' \emph{IEEE Transactions on Power Systems}, vol.~34, no.~6, pp. 4910--4920, Nov. 2019.

\bibitem{dobbe_linear_2020}
R.~Dobbe, W.~van Westering, S.~Liu, D.~Arnold, D.~Callaway, and C.~Tomlin, ``Linear single- and three-phase voltage forecasting and bayesian state estimation with limited sensing,'' \emph{IEEE Transactions on Power Systems}, vol.~35, no.~3, pp. 1674--1683, May 2020.

\bibitem{zhang_robust_2024}
J.~Zhang, Y.~Liu, A.~Rouhani, J.~Zhao, G.~J. Warchol, and K.~Scott, ``A robust data-driven gaussian process regression-enabled distribution system state estimation using partial node data,'' in \emph{2024 IEEE Power \& Energy Society General Meeting (PESGM)}, Jul. 2024, pp. 1--5, iSSN: 1944-9933.

\bibitem{Donti2020}
P.~L. Donti, Y.~Liu, A.~J. Schmitt, A.~Bernstein, R.~Yang, and Y.~Zhang, ``Matrix completion for low-observability voltage estimation,'' \emph{IEEE Transactions on Smart Grid}, vol.~11, no.~3, pp. 2520--2530, 2020.

\bibitem{liu_data_2025}
Y.~Liu, A.~Rouhani, J.~Zhao, J.~Zhang, G.~J. Warchol, and K.~Scott, ``Data enhanced robust distribution system state estimation for realistic system with parameter corrections,'' \emph{IEEE Transactions on Power Systems}, pp. 1--12, 2025.

\bibitem{zeraati_novel_2024}
M.~Zeraati, M.~R. Sheibani, F.~Jabari, and E.~Heydarian-Forushani, ``A novel state estimation method for distribution networks with low observability based on linear {AC} optimal power flow model,'' \emph{Electric Power Systems Research}, vol. 228, p. 110085, Mar. 2024.

\bibitem{karimi_joint_2021}
H.~S. Karimi and B.~Natarajan, ``Joint topology identification and state estimation in unobservable distribution grids,'' \emph{IEEE Transactions on Smart Grid}, vol.~12, no.~6, pp. 5299--5309, Nov. 2021.

\bibitem{pau_wls-based_2024}
M.~Pau and P.~A. Pegoraro, ``{WLS}-based state estimation for unobservable distribution grids through allocation factors evaluation,'' \emph{IEEE Transactions on Instrumentation and Measurement}, vol.~73, pp. 1--13, 2024.

\bibitem{zhang2023robust}
X.~Zhang, W.~Yan, M.~Huo, and H.~Li, ``Robust interval state estimation for distribution systems considering pseudo-measurement interval prediction,'' \emph{Journal of Modern Power Systems and Clean Energy}, vol.~12, no.~1, pp. 179--188, 2023.

\bibitem{vanin2023exact}
M.~Vanin, T.~Van~Acker, R.~D’hulst, and D.~Van~Hertem, ``Exact modeling of non-gaussian measurement uncertainty in distribution system state estimation,'' \emph{IEEE Transactions on Instrumentation and Measurement}, vol.~72, pp. 1--11, 2023.

\bibitem{abur2004}
A.~Abur and A.~Gómez-Expósito, \emph{Power system state estimation: theory and implementation}, ser. Power engineering.\hskip 1em plus 0.5em minus 0.4em\relax New York, NY: Marcel Dekker, 2004.

\bibitem{wang2019}
G.~Wang, G.~B. Giannakis, and J.~Chen, ``Robust and scalable power system state estimation via composite optimization,'' \emph{IEEE Transactions on Smart Grid}, vol.~10, no.~6, pp. 6137--6147, Nov. 2019.

\bibitem{caro2020state}
E.~Caro and A.~Hernández, ``State estimation in power systems based on a mathematical programming approach,'' in \emph{Advances in Electric Power and Energy}.\hskip 1em plus 0.5em minus 0.4em\relax John Wiley \& Sons, Ltd, 2020, pp. 23--57.

\bibitem{delaVarga2025}
J.~G. De~la Varga, S.~Pineda, J.~M. Morales, and A.~Porras, ``Learning-based state estimation in distribution systems with limited real-time measurements,'' \emph{Electric Power Systems Research}, vol. 241, p. 111268, 2025.

\bibitem{gomez-exposito2015}
A.~Gómez-Expósito, C.~Gómez-Quiles, and I.~Džafić, ``State estimation in two time scales for smart distribution systems,'' \emph{IEEE Transactions on Smart Grid}, vol.~6, no.~1, pp. 421--430, Jan. 2015.

\bibitem{esteban2022distributionally}
A.~Esteban-Pérez and J.~M. Morales, ``Distributionally robust stochastic programs with side information based on trimmings,'' \emph{Mathematical Programming}, vol. 195, no.~1, pp. 1069--1105, 2022.

\bibitem{bobo2021case123}
L.~Bobo, A.~Venzke, and S.~Chatzivasileiadis, ``Second-order cone relaxations of the optimal power flow for active distribution grids: Comparison of methods,'' \emph{International Journal of Electrical Power \& Energy Systems}, vol. 127, p. 106625, 2021.

\bibitem{case30asADN_1}
H.~H. Zeineldin, Y.~A.-R.~I. Mohamed, V.~Khadkikar, and V.~R. Pandi, ``A protection coordination index for evaluating distributed generation impacts on protection for meshed distribution systems,'' \emph{IEEE Transactions on Smart Grid}, vol.~4, no.~3, pp. 1523--1532, 2013.

\bibitem{case30asADN_2}
H.~Shad, M.~Gandomkar, and J.~Nikoukar, ``Analysis the securable operation and protection coordination indices on multi-objective sitting and sizing of synchronous distributed generations on distribution networks,'' \emph{IET Generation, Transmission \& Distribution}, vol.~17, no.~9, pp. 2169--2181, 2023.

\bibitem{pandapower2018}
L.~Thurner, A.~Scheidler, F.~Schäfer, J.-H. Menke, J.~Dollichon, F.~Meier, S.~Meinecke, and M.~Braun, ``Pandapower—an open-source python tool for convenient modeling, analysis, and optimization of electric power systems,'' \emph{IEEE Transactions on Power Systems}, vol.~33, no.~6, pp. 6510--6521, 2018.

\bibitem{pecanstreet}
\BIBentryALTinterwordspacing
``\BIBforeignlanguage{en-US}{Pecan {Street} {Dataport}}.'' [Online]. Available: \url{https://www.pecanstreet.org/dataport/}
\BIBentrySTDinterwordspacing

\bibitem{noauthor_groupoasysrobustcontextualse_2025}
\BIBentryALTinterwordspacing
{OASYS Research Group}, ``{RobustContextualSE},'' Feb. 2026. [Online]. Available: \url{https://github.com/groupoasys/RobustContextualSE}
\BIBentrySTDinterwordspacing

\bibitem{massignan2022bayesian}
J.~A.~D. Massignan, J.~B.~A. London, M.~Bessani, C.~D. Maciel, R.~Z. Fannucchi, and V.~Miranda, ``Bayesian inference approach for information fusion in distribution system state estimation,'' \emph{IEEE Transactions on Smart Grid}, vol.~13, no.~1, pp. 526--540, 2022.

\bibitem{pyomo}
W.~E. Hart, J.~P. Watson, and D.~L. Woodruff, ``Pyomo: Modeling and solving mathematical programs in python,'' \emph{Mathematical Programming Computation}, vol.~3, no.~3, pp. 219--260, 2011.

\bibitem{knitro}
R.~H. Byrd, J.~Nocedal, and R.~A. Waltz, ``K-{NITRO}: An integrated package for nonlinear optimization,'' \emph{Large-scale nonlinear optimization}, pp. 35--59, 2006.

\bibitem{Kfouri2026bddi}
R.~Kfouri, R.~A. Jabr, and I.~Džafić, ``Bad data detection and identification based on graph neural network for power system state estimation,'' \emph{Journal of Modern Power Systems and Clean Energy}, vol.~14, no.~2, pp. 760--772, 2026.

\end{thebibliography}

\begin{IEEEbiography}[{\includegraphics[width=1in,height=1.25in,clip,keepaspectratio]{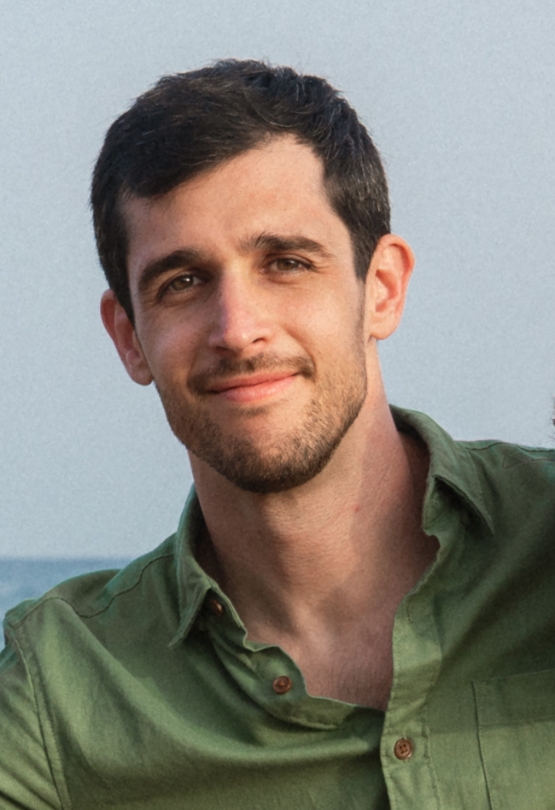}}]{Jos\'e G\'omez de la Varga } received the B.S. and M.S. degrees in industrial engineering from the University of Málaga, Spain, in 2016 and 2018, respectively, where he is currently pursuing the Ph.D. degree in electrical engineering with the OASYS Research Group. His research interests include distribution system state estimation, real-time observability, robust optimization, contextual stochastic optimization applied to active distribution networks, and physics-informed machine learning. 
\end{IEEEbiography}

\begin{IEEEbiography}
[{\includegraphics[width=1in,height=1.25in,clip,keepaspectratio]{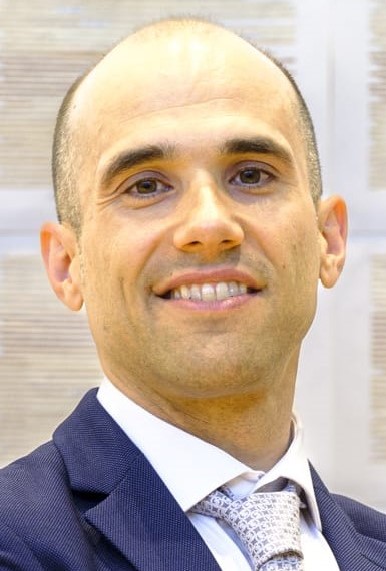}}]
{Juan Miguel Morales } received the B.Sc. and M.Sc. degrees in industrial engineering from the Univeristy of M\'alaga, M\'alaga, Spain, in 2006, and the Ph.D. degree in electrical engineering from the University of Castilla-La Mancha, Ciudad Real, Spain, in 2010.
He is currently a Full Professor with the Department of Mathematical Analysis, Statistics and Operations Research, and Applied Mathematics, University of M\'alaga. His research interests include fields of data analytics and optimization, decision-making under uncertainty, smart grids, power systems economics, operations and planning, and electricity markets.
\end{IEEEbiography}

\begin{IEEEbiography}
[{\includegraphics[width=1in,height=1.25in,clip,keepaspectratio]{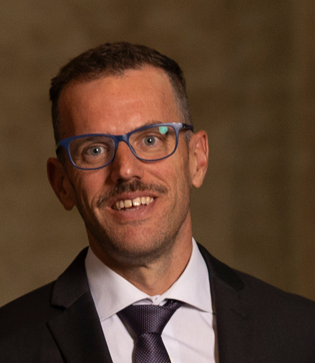}}]
{Salvador Pineda } received the B.Sc. and M.Sc. degrees in industrial engineering from the Univeristy of M\'alaga, M\'alaga, Spain, in 2006, and the Ph.D. degree in electrical engineering from the University of Castilla-La Mancha, Ciudad Real, Spain, in 2011.
He is currently a Full Professor with the Department of Electrical Engineering, University of M\'alaga. His research interests include power system operation and planning, decision-making under uncertainty, bilevel programming, machine learning, and statistics.
\end{IEEEbiography}

\end{document}